# Evaluation of integrals by differentiation with respect to a parameter

**Khristo N. Boyadzhiev**

Department of Mathematics and Statistics,
Ohio Northern University, Ada, OH 45810, USA
k-boyadzhiev@onu.edu

**Abstract**. We review a special technique for evaluating challenging integrals by providing a number of examples. Many of our examples prove integrals from the popular table of Gradshteyn and Ryzhik.

## 1 Introduction and three examples

There are various methods for evaluating integrals: substitution, integration by parts, partial fractions, using the residue theorem, or Cauchy's integral formula, etc. A beautiful special technique is the differentiation with respect to a parameter inside the integral. We review this technique here by providing numerous examples, many of which prove entries from the popular handbook of Gradshteyn and Ryzhik [6]. In our examples we focus on the formal manipulation. Several theorems justifying the legitimacy of the work are listed at the end of the paper. Applying the theorems in every particular case is left to the reader.

We hope integral lovers will appreciate this review and many will use it as a helpful reference. The examples and techniques are accessible to advanced calculus students and can be applied in various projects. Many more integrals from [6] can be proved by using the same approach.

We also want to mention that for many of the presented examples, solving the integral by differentiation with respect to a parameter is possibly the best solution.

Our main reference is the excellent book of Fikhtengolts [5] which is the source of several examples. Some integrals solved by this technique can be found in [1] and [2]. The method is presented in various publications, for instance, [4], [7], [9], [10], [12], and [13].

Below we evaluate three very different integrals in order to demonstrate the wide scope of the method. In section 2 we present a collection of eighteen more or less typical cases. In section 3 we show how differential equations can be involved very effectively. In section 4 we demonstrate a more sophisticated technique, where the parameter appears also in the integral limits.

**Example 1.1**

We start with a very simple example. It is easy to show that the popular integral

(1.1) $$\int_0^\infty \frac{\sin x}{x} dx$$

is convergent. In order to evaluate this integral we introduce the function

(1.2) $$F(\lambda) = \int_0^\infty e^{-\lambda x} \frac{\sin x}{x} dx, \ \lambda > 0$$

and differentiate this function to get

$$F'(\lambda) = -\int_0^\infty e^{-\lambda x} \sin x \, dx = \frac{-1}{1+\lambda^2}$$

(Laplace transform of the sine function). Integrating back we find

$$F(\lambda) = -\arctan(\lambda) + C \ .$$

Setting $\lambda \to \infty$ yields the equation

$$0 = -\frac{\pi}{2} + C, \text{ i.e. } C = \frac{\pi}{2},$$

and therefore,

(1.4) $$F(\lambda) = \int_0^\infty e^{-\lambda x} \frac{\sin x}{x} dx = \frac{\pi}{2} - \arctan \lambda$$

(cf. entry 3.9411 in [6]). Taking limits for $\lambda \to 0$ we find

(1.5) $$\int_0^\infty \frac{\sin x}{x} dx = \frac{\pi}{2}.$$

**Example 1.2**

The 66 Annual William Lowell Putnam Mathematical Competition (2005) included the integral (A5)

(1.6) $$\int_0^1 \frac{\ln(1+x)}{1+x^2} dx$$

with a solution published in [11]. This is entry 4.291(8) in [6]. We shall give a different solution by introducing a parameter. Consider the function

(1.7) $$F(\lambda) = \int_0^1 \frac{\ln(1+\lambda x)}{1+x^2} dx$$

defined for $\lambda \geq 0$. Differentiating this function we get

$$F'(\lambda) = \int_0^1 \frac{x}{(1+\lambda x)(1+x^2)} dx.$$

This integral is easy to evaluate by splitting the integrand in partial fractions. The result is

$$F'(\lambda) = -\frac{\ln(1+\lambda)}{1+\lambda^2} + \frac{1}{2} \ln 2 \frac{1}{1+\lambda^2} + \frac{\pi}{4} \frac{\lambda}{1+\lambda^2}.$$

Integrating we find

(1.8) $$F(\lambda) = -\int_0^\lambda \frac{\ln(1+x)}{1+x^2} dx + \frac{\ln 2}{2} \arctan \lambda + \frac{\pi}{8} \ln(1+\lambda^2),$$

and setting $\lambda = 1$ we arrive at the equation

$$2F(1) = \frac{\pi}{4} \ln 2.$$

That is,

(1.9) $$\int_0^1 \frac{\ln(1+x)}{1+x^2} dx = \frac{\pi}{8} \ln 2.$$

As we shall see later, many integrals containing logarithms and inverse trigonometric functions can be evaluated by this method.

It is good to notice that integrating (1.9) by parts we find

$$\int_0^1 \frac{\ln(1+x)}{1+x^2} dx = \ln(1+x) \arctan x \Big|_0^1 - \int_0^1 \frac{\arctan x}{1+x} dx$$

and therefore, we have also

$$\int_0^1 \frac{\arctan x}{1+x} dx = \frac{\pi}{8} \ln 2.$$

**Example 1.3**

This is Problem 1997 from the Mathematics Magazine 89(3), 2016, p. 223. Evaluate

(1.10) $$\int_0^\infty \left(\frac{1-e^{-x}}{x}\right)^2 dx.$$

Solution. We show that for every $\lambda > 0$

$$(1.11) \qquad F(\lambda) \equiv \int_0^\infty \left( \frac{1 - e^{-\lambda x}}{x} \right)^2 dx = \lambda \ln 4 \; .$$

Indeed, differentiating this function with respect to $\lambda$ (which is legitimate, as the integral is uniformly convergent on every interval $0 < a < \lambda < b$) we find

$$F'(\lambda) = 2 \int_0^\infty \left( \frac{1 - e^{-\lambda x}}{x} \right) e^{-\lambda x} dx = 2 \int_0^\infty \frac{e^{-\lambda x} - e^{-2\lambda x}}{x} dx = 2 \ln 2$$

by using Frullani's formula for the last equality (see below).

We conclude that $F(\lambda)$ is a linear function and since $F(0) = 0$ we can write $F(\lambda) = \lambda \ln 4$. With $\lambda = 1$ we find $F(1) = \ln 4$.

Frullani's formula says that for appropriate functions $f(x)$ we have

$$(1.12) \qquad \int_0^\infty \frac{f(ax) - f(bx)}{x} dx = [f(0) - f(\infty)] \ln \frac{b}{a} \; .$$

## 2. General examples

**Example 2.1**

We start this section with a simple and popular example Consider the integral

$$(2.1) \qquad J(\alpha) = \int_0^1 \frac{x^\alpha - 1}{\ln x} dx$$

for $\alpha \geq 0$. Note that the integrand is a continuous function on $[0,1]$ when it is defined as zero at $x = 0$ and as $\alpha$ (its limit value) at $x = 1$. Since $\frac{d}{d\alpha} x^\alpha = x^\alpha \ln x$, differentiation with respect to $\alpha$ yields

$$J'(\alpha) = \int_0^1 x^\alpha dx = \frac{1}{1 + \alpha}$$

and $J(\alpha) = \ln(1+\alpha) + C$. Since $J(0) = 0$ we find $C = 0$ and finally,

(2.2) $$\int_0^1 \frac{x^\alpha - 1}{\ln x} dx = \ln(1+\alpha).$$

The similar integral

(2.3) $$\int_0^1 \frac{x^\alpha - x^\beta}{\ln x} dx$$

where $\alpha, \beta \geq 0$ can be reduced to (2.2) by writing $x^\alpha - x^\beta = x^\alpha - 1 - (x^\beta - 1)$. Thus we have

(2.4) $$\int_0^1 \frac{x^\alpha - x^\beta}{\ln x} dx = \ln \frac{1+\alpha}{1+\beta}.$$

Another way to approach (2.1) is to use the substitution $x = e^{-t}$ which transforms it to the Frullani integral

(2.5) $$\int_0^\infty \frac{e^{-(\alpha+1)t} - e^{-t}}{t} dt$$

see (1.12).

**Example 2.2**

We evaluate here the improper integral

(2.6) $$J(\lambda) = \int_0^1 \frac{\arctan \lambda x}{x\sqrt{1-x^2}} dx.$$

Differentiation yields

$$J'(\lambda) = \int_0^1 \frac{dx}{(1+\lambda^2 x^2)\sqrt{1-x^2}}$$

and with the substitution $x = \cos\theta$ this transforms into

$$J'(\lambda) = \int_0^{\pi/2} \frac{1}{1+\lambda^2 \cos^2 \theta} d\theta = \int_0^{\pi/2} \frac{d\theta}{(1+\tan^2 \theta + \lambda^2)\cos^2 \theta}$$

$$= \int_0^{\pi/2} \frac{d\tan\theta}{1+\lambda^2 + \tan^2 \theta} = \frac{1}{\sqrt{1+\lambda^2}} \arctan \frac{\tan\theta}{\sqrt{1+\lambda^2}} \bigg|_0^{\frac{\pi}{2}} = \frac{\pi}{2\sqrt{1+\lambda^2}}.$$

Therefore,

(2.7) $$J(\lambda) = \frac{\pi}{2} \ln\left(\lambda + \sqrt{1+\lambda^2}\right),$$

since $J(0) = 0$. In particular, with $\lambda = 1$,

(2.8) $$\int_0^1 \frac{\arctan x}{x\sqrt{1-x^2}} dx = \frac{\pi}{2} \ln(1+\sqrt{2}).$$

This integral is entry 4.531(12) in [6]. Note that the similar integral

(2.9) $$J(\lambda) = \int_0^1 \frac{\arctan \lambda x}{\sqrt{1-x^2}} dx$$

cannot be evaluated in the same manner. The derivative here becomes

$$J'(\lambda) = \int_0^1 \frac{x\,dx}{(1+\lambda^2 x^2)\sqrt{1-x^2}} = \int_0^{\pi/2} \frac{\cos\theta}{1+\lambda^2 \cos^2 \theta} d\theta = \frac{1}{2\lambda\sqrt{1+\lambda^2}} \ln \frac{\sqrt{1+\lambda^2}+\lambda}{\sqrt{1+\lambda^2}-\lambda}$$

which is not easy to integrate. The integral (2.9) will be evaluated later in section 4 by a more sophisticated method.

**Example 2.3**

Now consider

(2.10) $$J(\lambda) = \int_0^\infty \frac{\ln(1+\lambda^2 x^2)}{1+x^2} dx$$

with

$$J'(\lambda) = \int_0^\infty \frac{2\lambda x^2}{(1+\lambda^2 x^2)(1+x^2)} dx$$

$$= \frac{2\lambda}{1-\lambda^2} \int_0^\infty \left( \frac{1}{1+\lambda^2 x^2} - \frac{1}{1+x^2} \right) dx = \frac{2\lambda}{1-\lambda^2} \left( \frac{\pi}{2\lambda} - \frac{\pi}{2} \right) = \frac{\pi}{1+\lambda}$$

under the restriction $\lambda^2 \neq 1$. This way

(2.11) $\qquad J(\lambda) = \pi \ln(1+\lambda)$.

We needed $\lambda^2 \neq 1$ for the evaluation of $J'(\lambda)$, but this restriction can later be dropped. For equation (2.11) we only need $\lambda > -1$. In particular, for $\lambda = 1$,

$$\int_0^\infty \frac{\ln(1+x^2)}{1+x^2} dx = \pi \ln 2 .$$

The similar integral

$$\int_0^1 \frac{\ln(1+x^2)}{1+x} dx = \frac{3}{4}(\ln 2)^2 - \frac{\pi^2}{48}$$

is evaluated by the same method in [3]. In that article one can find also the evaluation

$$\int_0^1 \frac{\ln(1+x^2)}{1+x^2} dx = \frac{\pi}{2} \ln 2 - G ,$$

where $G$ is Catalan's constant (see the remark at the end of Section 4).

**Example 2.4**

We shall evaluate here two more integrals with arctangents. The first one is 4.535(7) from [6],

(2.12) $\qquad G(\lambda) = \int_0^\infty \frac{\arctan \lambda x}{x(1+x^2)} dx$ .

For all $\lambda > -1, \lambda \neq 1$ we compute

$$G'(\lambda) = \int_0^\infty \frac{dx}{(1+\lambda^2 x^2)(1+x^2)} = \frac{1}{1-\lambda^2} \int_0^\infty \left( \frac{1}{1+x^2} - \frac{\lambda^2}{1+\lambda^2 x^2} \right) dx$$

$$= \frac{1}{1-\lambda^2} \left( \frac{\pi}{2} - \frac{\pi \lambda}{2} \right) = \frac{\pi}{2(1+\lambda)},$$

and hence (dropping the restriction $\lambda \neq 1$)

(2.13) $$G(\lambda) = \frac{\pi}{2} \ln(1+\lambda).$$

and for $\lambda = 1$

(2.14) $$\int_0^\infty \frac{\arctan x}{x(1+x^2)} dx = \frac{\pi}{2} \ln 2.$$

Comparing this to (2.11) we conclude that for all $\lambda > -1$

(2.15) $$\int_0^\infty \frac{\ln(1+\lambda^2 x^2)}{1+x^2} dx = 2 \int_0^\infty \frac{\arctan \lambda x}{x(1+x^2)} dx = \pi \ln(1+\lambda).$$

**Example 2.5**

Related to (2.12) is the following integral

(2.16) $$G(\lambda, \mu) = \int_0^\infty \frac{\arctan(\lambda x) \arctan(\mu x)}{x^2} dx$$

for $\lambda, \mu > 0$. Using the evaluation (2.13) we write

$$G_\lambda(\lambda, \mu) = \int_0^\infty \frac{\arctan(\mu x)}{x(1+\lambda^2 x^2)} dx = \frac{\pi}{2} \ln\left(1 + \frac{\lambda}{\mu}\right)$$

and integrating this logarithm by parts with respect to $\lambda$ we find

$$G(\lambda, \mu) = \frac{\pi}{2} [(\lambda + \mu) \ln(\lambda + \mu) - \lambda \ln \lambda] + C(\mu).$$

Setting $\lambda \to 0$ yields $C(\mu) = -\dfrac{\pi}{2}\mu \ln \mu$. Finally,

(2.17) $$\int_0^\infty \frac{\arctan(\lambda x)\arctan(\mu x)}{x^2} dx = \frac{\pi}{2}[(\lambda + \mu)\ln(\lambda + \mu) - \lambda \ln \lambda - \mu \ln \mu].$$

**Example 2.6**

Consider the integral 3.943 from [6]

$$F(\lambda) = \int_0^\infty e^{-\beta x}\frac{1-\cos \lambda x}{x}dx,$$

where $\beta > 0$ is fixed. We have

$$F'(\lambda) = \int_0^\infty e^{-\beta x}\sin \lambda x\, dx = \frac{\lambda}{\lambda^2 + \beta^2}$$

and integrating back

$$F(\lambda) = \frac{1}{2}\ln(\lambda^2 + \beta^2) + C(\beta).$$

To compute $C(\beta)$ we set $\lambda = 0$ and this gives $C(\beta) = \dfrac{-1}{2}\ln \beta^2$. Therefore,

(2.18) $$\int_0^\infty e^{-\beta x}\frac{1-\cos \lambda x}{x}dx = \frac{1}{2}\ln\left(1 + \frac{\lambda^2}{\beta^2}\right).$$

**Example 2.7**

A "symmetrical" analog to the previous example is the integral

$$F(\lambda) = \int_0^\infty \frac{1 - e^{-\lambda x}}{x}\cos \beta x\, dx,$$

defined for $\lambda \geq 0$ and $\beta \neq 0$. The integral is divergent at infinity when $\beta = 0$. We have

$$F'(\lambda) = \int_0^\infty e^{-\lambda x} \cos \beta x \, dx = \frac{\lambda}{\lambda^2 + \beta^2}$$

and integrating

(2.19) $$F(\lambda) = \int_0^\infty \frac{1 - e^{-\lambda x}}{x} \cos \beta x \, dx = \frac{1}{2} \ln\left(1 + \frac{\lambda^2}{\beta^2}\right)$$

so that for any $\lambda \geq 0, \beta > 0$

$$\int_0^\infty \frac{1 - e^{-\lambda x}}{x} \cos \beta x \, dx = \int_0^\infty e^{-\beta x} \frac{1 - \cos \lambda x}{x} dx.$$

Note that the integral 3.951(3) from [6]

$$\int_0^\infty \frac{e^{-\lambda x} - e^{-\mu x}}{x} \cos \beta x \, dx$$

can be reduced to (2.19) by writing $e^{-\lambda x} - e^{-\mu x} = (e^{-\lambda x} - 1) + (1 - e^{-\mu x})$ and splitting it in two integrals. Thus

$$\int_0^\infty \frac{e^{-\lambda x} - e^{-\mu x}}{x} \cos \beta x \, dx = \frac{1}{2} \ln \frac{\mu^2 + \beta^2}{\lambda^2 + \beta^2}.$$

**Example 2.8**

Using the well-known *Gaussian integral*, also known as the *Euler-Poisson integral*,

(2.20) $$\int_0^\infty e^{-x^2} dx = \frac{\sqrt{\pi}}{2}$$

we can evaluate for every $\lambda \geq 0$ the integral

$$F(\lambda) = \int_0^\infty \frac{1 - e^{-\lambda x^2}}{x^2} dx .$$

Indeed, we have for $\lambda > 0$

$$F'(\lambda) = \int_0^\infty e^{-\lambda x^2}\,dx = \frac{1}{\sqrt{\lambda}}\int_0^\infty \exp(-(x\sqrt{\lambda})^2)\,dx\sqrt{\lambda} = \frac{\sqrt{\pi}}{2\sqrt{\lambda}},$$

so that

(2.21) $$F(\lambda) = \int_0^\infty \frac{1 - e^{-\lambda x^2}}{x^2}\,dx = \sqrt{\lambda\pi}.$$

**Example 2.9**

Sometimes we can use partial derivatives as in the following integral. Consider the function

(2.22) $$F(\lambda, \mu) = \int_0^\infty \frac{e^{-px}\cos qx - e^{-\lambda x}\cos \mu x}{x}\,dx$$

with four parameters. Here $\lambda > 0$, $\mu$ will be variables and $p > 0$, $q$ will be fixed. The partial derivatives are

$$F_\lambda(\lambda, \mu) = \int_0^\infty e^{-\lambda x}\cos \mu x\,dx = \frac{\lambda}{\lambda^2 + \mu^2},$$

$$F_\mu(\lambda, \mu) = \int_0^\infty e^{-\lambda x}\sin \mu x\,dx = \frac{\mu}{\lambda^2 + \mu^2}.$$

It is easy to restore the function from these derivatives

$$F(\lambda, \mu) = \frac{1}{2}\ln(\lambda^2 + \mu^2) + C(p, q),$$

where $C(p, q)$ is unknown. The integral (2.22) is zero when $\lambda = p$ and $\mu = q$, so from the last equation we find $C(p, q) = -\ln(p^2 + q^2)/2$. Therefore,

(2.23) $$\int_0^\infty \frac{e^{-px}\cos qx - e^{-\lambda x}\cos \mu x}{x}\,dx = \frac{1}{2}\ln\frac{\lambda^2 + \mu^2}{p^2 + q^2}.$$

**In all following examples containing "$e^{-\lambda x}$" we assume $\lambda > 0$.**

**Example 2.10**

Now consider

$$J(\lambda) = \int_0^\infty e^{-\lambda x} \frac{\sin(ax)\sin(bx)}{x} dx$$

where $a > b > 0$ are constants. Clearly,

$$J'(\lambda) = -\int_0^\infty e^{-\lambda x} \sin(ax)\sin(bx) dx$$

$$= \frac{1}{2}\left\{\int_0^\infty e^{-\lambda x} \cos(a+b)x\, dx - \int_0^\infty e^{-\lambda x} \cos(a-b)x\, dx\right\}$$

$$= \frac{1}{2}\left\{\frac{\lambda}{\lambda^2 + (a+b)^2} - \frac{\lambda}{\lambda^2 + (a-b)^2}\right\}.$$

Integrating with respect to $\lambda$ and evaluating the constant of integration with $\lambda \to \infty$ we find

(2.24) $$J(\lambda) = \int_0^\infty e^{-\lambda x} \frac{\sin(ax)\sin(bx)}{x} dx = \frac{1}{4}\ln\frac{\lambda^2 + (a+b)^2}{\lambda^2 + (a-b)^2}.$$

This is entry 3.947(1) in [6].

**Example 2.11**

Using the previous example we can evaluate also entry 3.947(2) in [6].

$$G(\lambda) = \int_0^\infty e^{-\lambda x} \frac{\sin(ax)\sin(bx)}{x^2} dx$$

where again $a > b > 0$. We have from above

$$G'(\lambda) = -J(\lambda) = \frac{-1}{4}\ln\frac{\lambda^2 + (a+b)^2}{\lambda^2 + (a-b)^2} = \frac{1}{4}\ln\frac{\lambda^2 + (a-b)^2}{\lambda^2 + (a+b)^2}$$

and integrating by parts,

$$G(\lambda) = \frac{\lambda}{4} \ln \frac{\lambda^2 + (a-b)^2}{\lambda^2 + (a+b)^2} - \frac{1}{4}\int\left(\frac{\lambda^2}{\lambda^2+(a-b)^2} - \frac{\lambda^2}{\lambda^2+(a+b)^2}\right)d\lambda.$$

With simple algebra we find

$$\frac{\lambda^2}{\lambda^2+(a-b)^2} - \frac{\lambda^2}{\lambda^2+(a+b)^2} = \frac{(a+b)^2}{\lambda^2+(a+b)^2} - \frac{(a-b)^2}{\lambda^2+(a-b)^2}$$

and the integration becomes easy. The result is

$$G(\lambda) = \frac{\lambda}{4}\ln\frac{\lambda^2+(a-b)^2}{\lambda^2+(a+b)^2} + \frac{a-b}{2}\arctan\frac{\lambda}{a-b} - \frac{a+b}{2}\arctan\frac{\lambda}{a+b} + \frac{\pi b}{2}$$

(the constant of integration is found by setting $\lambda \to \infty$).

This answer is simpler than the one given in [6]. With $\lambda = 0$ we prove also 3.741(3) from [6]

$$\int_0^\infty \frac{\sin(ax)\sin(bx)}{x^2}dx = \frac{\pi b}{2} \quad (a \geq b > 0).$$

**Example 2.12**

Similar to (2.24) is the integral

$$G(\lambda) = \int_0^\infty e^{-\lambda x}\frac{\sin(ax)\cos(bx)}{x}dx$$

(this is 3.947(3) in [6]). Suppose $a > b > 0$. Then

$$G'(\lambda) = -\int_0^\infty e^{-\lambda x}\sin(ax)\cos(bx)\,dx$$

$$= \frac{-1}{2}\left\{\int_0^\infty e^{-\lambda x}\sin(a+b)x\,dx + \int_0^\infty e^{-\lambda x}\sin(a-b)x\,dx\right\}$$

$$= \frac{-1}{2}\left\{\frac{a+b}{\lambda^2+(a+b)^2} + \frac{a-b}{\lambda^2+(a-b)^2}\right\},$$

and after integration with respect to $\lambda$,

$$G(\lambda) = \frac{\pi}{2} - \frac{1}{2}\left(\arctan\frac{\lambda}{a+b} + \arctan\frac{\lambda}{a-b}\right),$$

where the constant of integration $\pi/2$ is found by letting $\lambda \to \infty$.

Setting $b \to a$ we find also

$$\int_0^\infty e^{-\lambda x} \frac{\sin(ax)\cos(ax)}{x} dx = \frac{\pi}{4} - \frac{1}{2}\arctan\frac{\lambda}{2a}.$$

Using the identity $2\sin(ax)\cos(ax) = \sin(2ax)$ this integral can be reduced to (1.4).

**Example 2.13**

(2.25) $$F(\lambda) = \int_0^\infty e^{-\lambda x} \frac{\cos(ax) - \cos(bx)}{x^2} dx$$

(entry 3.948(3) in [6]). Differentiating we find

$$F'(\lambda) = \int_0^\infty e^{-\lambda x} \frac{\cos(bx) - \cos(ax)}{x} dx = \frac{1}{2}\ln\frac{\lambda^2 + a^2}{\lambda^2 + b^2}$$

in view of (2.18), as $\cos(bx) - \cos(ax) = \cos(bx) - 1 + 1 - \cos(ax)$. Integration by parts yields

(2.26) $$F(\lambda) = \frac{\lambda}{2}\ln\frac{\lambda^2 + a^2}{\lambda^2 + b^2} + b\arctan\frac{b}{\lambda} - a\arctan\frac{a}{\lambda}$$

(again the constant of integration is found by setting $\lambda \to \infty$).

Various integrals with similar structure can be evaluated by this method or by reducing to those already evaluated here. For example, entry 3.948 (4) from [6]

$$A(\lambda) = \int_0^\infty e^{-\lambda x} \frac{\sin^2(ax) - \sin^2(bx)}{x^2} dx,$$

can be reduced to (2.25) by using the identity $2\sin^2 x = 1 - \cos 2x$. The results is

$$A(\lambda) = \frac{\lambda}{4} \ln \frac{\lambda^2 + 4b^2}{\lambda^2 + 4a^2} + a \arctan \frac{2a}{\lambda} - b \arctan \frac{2b}{\lambda} \ .$$

Now we shall evaluate several integrals involving logarithms of trigonometric functions.

**Example 2.14**

Consider the integral

(2.27) $$J(\alpha) = \int_0^{\frac{\pi}{2}} \ln(\alpha^2 - \cos^2 \theta) d\theta$$

for $\alpha > 1$. Differentiation with respect to $\alpha$ yields

$$J'(\alpha) = 2\alpha \int_0^{\frac{\pi}{2}} \frac{d\theta}{\alpha^2 - \cos^2 \theta}$$

and then the substitution $x = \tan \theta$ turns this into

$$J'(\alpha) = 2\alpha \int_0^\infty \frac{dx}{\alpha^2 - 1 + \alpha^2 x^2} = \frac{2}{\sqrt{\alpha^2 - 1}} \arctan \frac{\alpha x}{\sqrt{\alpha^2 - 1}} \bigg|_0^\infty = \frac{\pi}{\sqrt{\alpha^2 - 1}} \ .$$

Therefore,

$$J(\alpha) = \int_0^{\frac{\pi}{2}} \ln(\alpha^2 - \cos^2 \theta) d\theta = \pi \ln\left(\alpha + \sqrt{\alpha^2 - 1}\right) + C \ .$$

In order to evaluate the constant of integration we write this equation in the form (factoring out $\alpha^2$ in the left hand side and $\alpha$ in the right hand side)

$$\pi \ln \alpha + \int_0^{\frac{\pi}{2}} \ln\left(1 - \frac{\cos^2 \theta}{\alpha^2}\right) d\theta = \pi \ln \alpha + \pi \ln\left(1 + \sqrt{1 - \frac{1}{\alpha^2}}\right) + C \ .$$

Removing $\pi \ln \alpha$ from both sides and setting $\alpha \to \infty$ we compute $C = -\pi \ln 2$. As a result, two integrals are evaluated. For the second one we set $\beta = 1/\alpha$ in (2.27)

(2.28) $$\int_0^{\frac{\pi}{2}} \ln(\alpha^2 - \cos^2\theta)d\theta = \pi \ln\frac{\alpha + \sqrt{\alpha^2-1}}{2} \quad (\alpha > 1)$$

(2.29) $$\int_0^{\frac{\pi}{2}} \ln(1 - \beta^2 \cos^2\theta)d\theta = \pi \ln\frac{1 + \sqrt{1-\beta^2}}{2} \quad (0 \le \beta \le 1).$$

In particular, with $\beta = 1$ in (2.29) we obtain the important log-sine integral

(2.30) $$\int_0^{\frac{\pi}{2}} \ln(\sin\theta)d\theta = -\frac{\pi}{2}\ln 2.$$

**Example 2.15**

In the same way we can prove that

(2.31) $$\int_0^{\pi/2} \ln(1 + \alpha \sin^2\theta)d\theta = \pi \ln\frac{1 + \sqrt{1+\alpha}}{2}$$

for any $\alpha > -1$. Calling this integral $F(\alpha)$ and differentiating we find

$$F'(\alpha) = \int_0^{\pi/2} \frac{\sin^2\theta}{1 + \alpha \sin^2\theta}d\theta.$$

Now we divide top and bottom of the integrand by $\cos^2\theta$ and then use the substitution $x = \tan\theta$

$$F'(\alpha) = \int_0^{\infty} \frac{x^2}{(1+x^2)(1+(1+\alpha)x^2)}dx.$$

Assuming for the moment that $\alpha \ne 0$ and using partial fraction we write

$$F'(\alpha) = \frac{1}{\alpha}\int_0^{\infty}\left(\frac{1}{1+x^2} - \frac{1}{1+(1+\alpha)x^2}\right)dx$$

$$= \frac{1}{\alpha}\left(\arctan x - \frac{1}{\sqrt{1+\alpha}}\arctan(x\sqrt{1+\alpha})\right)\Big|_0^{\infty} = \frac{\pi}{2\alpha}\left(1 - \frac{1}{\sqrt{1+\alpha}}\right) = \frac{\pi}{2}\frac{\sqrt{1+\alpha}-1}{\alpha\sqrt{1+\alpha}}.$$

Simple algebra shows that

$$F'(\alpha) = \frac{\pi}{2(1+\sqrt{1+\alpha})\sqrt{1+\alpha}}$$

which is exactly the derivative of $\pi \ln(1+\sqrt{1+\alpha})$. Thus $F(\alpha) = \pi \ln(1+\sqrt{1+\alpha}) + C$. Setting $\alpha = 0$ we find $C = -\pi \ln 2$ and (2.31) is proved.

**Example 2.16**

Let $|\alpha| < 1$. Now we prove the interesting integral, entry 4.397 (3) in [6]

(2.32) $$F(\alpha) \equiv \int_0^\pi \frac{\ln(1+\alpha \cos \theta)}{\cos \theta} d\theta = \pi \arcsin \alpha.$$

Assuming that the value of the integrand at $\theta = \pi/2$ is $\alpha$, the integrand becomes a continuous function on $[0, \pi]$. Then

$$F'(\alpha) \equiv \int_0^\pi \frac{1}{1+\alpha \cos \theta} d\theta$$

which is easily solved with the substitution $\tan \frac{\theta}{2} = t$, so that $\cos \theta = \frac{1-t^2}{1+t^2}$, $d\theta = \frac{2dt}{1+t^2}$. Thus

$$F'(\alpha) = 2\int_0^\infty \frac{dt}{1+\alpha+(1-\alpha)t^2} = \frac{2}{1+\alpha}\int_0^\infty \frac{dt}{1+\frac{1-\alpha}{1+\alpha}t^2}$$

$$= \frac{2}{\sqrt{1-\alpha^2}} \arctan\left(t\sqrt{\frac{1-\alpha}{1+\alpha}}\right)\Bigg|_0^\infty = \frac{\pi}{\sqrt{1-\alpha^2}},$$

and (2.32) follows since both sides in this equation are zeros for $\alpha = 0$.

**Example 2.17**

Let again $|\alpha| < 1$. Using the results from the previous example we can prove

(2.33) $$\int_0^\pi \ln(1+\alpha\cos\theta)d\theta = \pi\ln\frac{1+\sqrt{1-\alpha^2}}{2}.$$

Indeed, let this integral be $F(\alpha)$. We have

$$F'(\alpha) = \int_0^\pi \frac{\cos\theta}{1+\alpha\cos\theta}d\theta = \frac{1}{\alpha}\int_0^\pi \frac{1+\alpha\cos\theta-1}{1+\alpha\cos\theta}d\theta = \frac{\pi}{\alpha} - \frac{1}{\alpha}\int_0^\pi \frac{d\theta}{1+\alpha\cos\theta}$$

$$= \frac{\pi}{\alpha} - \frac{\pi}{\alpha\sqrt{1-\alpha^2}}$$

(for the moment $\alpha \neq 0$). Integration is easy:

$$F(\alpha) = \pi\left(\ln\alpha + \int\frac{d\alpha^{-1}}{\sqrt{\alpha^{-2}-1}}\right) = \pi\left(\ln\alpha + \ln(\alpha^{-1}+\sqrt{\alpha^{-2}-1})\right) + C$$

$$= \pi\ln(1+\sqrt{1-\alpha^2}) + C.$$

Now we can drop the restriction $\alpha \neq 0$. With $\alpha = 0$ we find $C = -\pi\ln 2$ and (2.33) is proved.

**Example 2.18**

The last example in this section is a very interesting integral. It can be found, for example, in the book [12] on p. 143.

(2.34) $$F(\alpha) = \int_0^\pi \ln(1 - 2\alpha\cos x + \alpha^2)\,dx$$

We first assume $\alpha \neq 0$ and $\alpha \neq 1$. Then

$$F'(\alpha) = \int_0^\pi \frac{-2\cos x + 2\alpha}{1 - 2\alpha\cos x + \alpha^2}dx = \frac{1}{\alpha}\int_0^\pi\left(1 - \frac{1-\alpha^2}{1-2\alpha\cos x + \alpha^2}\right)dx$$

$$= \frac{\pi}{\alpha} - \frac{1-\alpha^2}{\alpha}\int_0^\pi \frac{1}{1-2\alpha\cos x + \alpha^2}dx.$$

The last integral can be evaluated by setting as before $\tan\dfrac{x}{2} = t$, with $\cos x = \dfrac{1-t^2}{1+t^2}$, $dx = \dfrac{2dt}{1+t^2}$.

Some simple work gives

$$F'(\alpha) = \dfrac{\pi}{\alpha} - \dfrac{2}{\alpha}\arctan\dfrac{1+\alpha}{1-\alpha}t\bigg|_0^\infty$$

and we find from here that $F'(\alpha) = 0$ when $|\alpha| < 1$ and $F'(\alpha) = \dfrac{2\pi}{\alpha}$ when $|\alpha| > 1$.

Thus $F(\alpha) = C_1$ when $|\alpha| < 1$ and $F(\alpha) = C_2 + \pi \ln \alpha^2$ for $|\alpha| > 1$. Since $F(0) = 0$ (as the definition (2.30) shows) we have

$$F(\alpha) = 0 \text{ for all } |\alpha| < 1.$$

Next, to determine the constant $C_2$ when $|\alpha| > 1$ we factor out $\alpha^2$ inside the logarithm in (2.34) and write

$$F(\alpha) = \int_0^\pi \ln\left(\alpha^2\left(\dfrac{1}{\alpha^2} - \dfrac{2\cos x}{\alpha} + 1\right)\right) dx = \pi \ln \alpha^2 + F\left(\dfrac{1}{\alpha}\right) = \pi \ln \alpha^2 + 0 = \pi \ln \alpha^2,$$

that is, $C_2 = 0$ and $F(\alpha) = \pi \ln \alpha^2$ for $|\alpha| > 1$. This also extends to $\alpha = \pm 1$, i.e.

$$F(\alpha) = \pi \ln \alpha^2 \text{ for } |\alpha| \geq 1.$$

It is good to mention here that the evaluation of this integral for the case $|\alpha| < 1$ can be done immediately by using the series representation for $0 \leq x \leq \pi$

(2.35) $$\ln(1 - 2\alpha \cos x + \alpha^2) = -2 \sum_{k=1}^\infty \dfrac{\alpha^k \cos kx}{k}.$$

The case $|\alpha| > 1$ can be reduced to this one by writing

$$\ln(1 - 2\alpha \cos x + \alpha^2) = \ln\left[\alpha^2\left(\dfrac{1}{\alpha^2} - 2\dfrac{1}{\alpha}\cos x + 1\right)\right] = 2\ln|\alpha| + \ln\left(\dfrac{1}{\alpha^2} - 2\dfrac{1}{\alpha}\cos x + 1\right).$$

The integral can be written in a symmetric form with $|\beta| \leq |\alpha|$ (cf. entry 2.6.36(14) in [8])

(2.36) $$\int_0^\pi \ln(\beta^2 - 2\alpha\beta \cos x + \alpha^2)\, dx = 2\pi \ln |\alpha|.$$

# 3. Using differential equations

**Example 3.1** Consider the integral

$$y(x) = \int_0^\infty e^{-t^2} \cos(2xt)\, dt.$$

Here

$$y'(x) = -\int_0^\infty 2t e^{-t^2} \sin(2xt)\, dt,$$

and integration by parts leads to the separable differential equation

$$y' = -2xy \quad \text{or} \quad \frac{dy}{dx} = -2xy$$

with general solution

$$y(x) = Ce^{-x^2}.$$

For $x = 0$ in the original integral we have $y(0) = \frac{\sqrt{\pi}}{2}$ according to (2.20). Therefore,

$$y(x) = \int_0^\infty e^{-t^2} \cos(2xt)\, dt = \frac{\sqrt{\pi}}{2} e^{-x^2}.$$

With a simple rescaling of the variable we can write this result as

$$\int_{-\infty}^\infty e^{-ax^2} \cos(xt)\, dx = \sqrt{\frac{\pi}{a}} e^{-t^2/4a} \quad (a > 0).$$

This last integral was used in the solution of Problem 1896 of the *Mathematics Magazine* (vol. 83, June 2013, 228-230).

**Example 3.2**

In this example we evaluate two interesting integrals (3.723, (2) and (3) in [6])

$$F(\lambda) = \int_0^\infty \frac{\cos \lambda x}{a^2 + x^2} dx, \text{ and } G(\lambda) = \int_0^\infty \frac{x \sin \lambda x}{a^2 + x^2} dx,$$

which can be viewed as Fourier cosine and sine transforms. We shall use a second order differential equation for $F(\lambda)$. First we have

(3.1) $\qquad F'(\lambda) = -G(\lambda)$

We cannot differentiate further, because $G'(\lambda)$ is divergent. Instead, we shall use a special trick, adding to both sides of (3.1) the number

$$\frac{\pi}{2} = \int_0^\infty \frac{\sin x}{x} dx$$

(see (1.5)). After a simple calculation

$$F'(\lambda) + \frac{\pi}{2} = a^2 \int_0^\infty \frac{\sin \lambda x}{x(a^2 + x^2)} dx.$$

Differentiating again we come to the second order differential equation

$$F'' = a^2 F$$

with general solution

$$F(\lambda) = Ae^{a\lambda} + Be^{-a\lambda},$$

where $A, B$ are arbitrary constants. Suppose $a > 0$ and $\lambda \geq 0$. Then $A = 0$, because $F(\lambda)$ is a bounded function when $\lambda \to \infty$. To find $B$ we set $\lambda = 0$ and use the fact that

$$B = F(0) = \int_0^\infty \frac{dx}{a^2+x^2} = \frac{1}{a}\arctan\frac{x}{a}\Big|_0^\infty = \frac{\pi}{2a}.$$

Finally,

(3.2) $$F(\lambda) = \int_0^\infty \frac{\cos\lambda x}{a^2+x^2}dx = \frac{\pi}{2a}e^{-a\lambda}$$

and from (3.1) we find also

(3.3) $$G(\lambda) = \int_0^\infty \frac{x\sin\lambda x}{a^2+x^2}dx = \frac{\pi}{2}e^{-a\lambda}.$$

This result can be used to evaluate some similar integrals. Integrating (3.2) with respect to $\lambda$ and adjusting the constant of integration we find entry 3.725 (1) [6]

$$\int_0^\infty \frac{\sin\lambda x}{x(a^2+x^2)}dx = \frac{\pi}{2a^2}(1-e^{-a\lambda}).$$

Differentiating this integral with respect to $a$ we prove also entry 3.735

$$\int_0^\infty \frac{\sin\lambda x}{x(a^2+x^2)^2}dx = \frac{\pi}{2a^4}(1-e^{-a\lambda}) - \frac{\lambda\pi}{4a^3}e^{-a\lambda}.$$

**Example 3.3**

We shall evaluate two Laplace integrals. For $s>0$ and $a>0$ consider

$$F(s) = \int_0^\infty \frac{e^{-st}}{a^2+t^2}dt, \text{ and } G(s) = \int_0^\infty \frac{te^{-st}}{a^2+t^2}dt.$$

Differentiating twice the first one we find

(3.4) $$F'(s) = -G(s),\ F'' = -G'(s)$$

and at the same time

$$-G'(s) = \int_0^\infty \frac{t^2 e^{-st}}{a^2+t^2} dt = \int_0^\infty \frac{(-a^2+a^2+t^2)e^{-st}}{a^2+t^2} dt = -a^2 F(s) + \frac{1}{s}$$

which leads to the second order differential equation

$$F'' + a^2 F = \frac{1}{s} .$$

This equation can be solved by variation of parameters. The solution is

$$F(s) = \frac{1}{a}[\text{ci}(as)\sin(as) - \text{si}(as)\cos(as)]$$

involving the special sine and cosine integrals

$$\text{si}(x) = -\int_x^\infty \frac{\sin t}{t} dt = \frac{-\pi}{2} + \int_0^x \frac{\sin t}{t} dt ,$$

$$\text{ci}(x) = -\int_x^\infty \frac{\cos t}{t} dt .$$

The choice of integral limits here is dictated by the initial conditions $F(\infty) = G(\infty) = 0$.

From (3.4) we find also

$$G(s) = -\text{ci}(as)\cos(as) - \text{si}(as)\sin(as) .$$

The integral $F(s)$ can be used to give an interesting extension of the integral (1.5). For any $a, b > 0$ we compute

$$\int_0^\infty \frac{\sin(ax)}{x+b} dx = \int_0^\infty \sin(ax)\left\{\int_0^\infty e^{-t(x+b)} dt\right\} dx = \int_0^\infty \left\{\int_0^\infty e^{-xt}\sin(ax) dx\right\} e^{-bt} dt$$

$$= a\int_0^\infty \frac{e^{-bt} dt}{t^2+a^2} = aF(b) = \text{ci}(ab)\sin(ab) - \text{si}(ab)\cos(ab) .$$

This is entry 3.772 (1) from [6]. Taking limits of both sides when $b \to 0$ yields (1.5). In the same way we prove entry 3.722(3)

$$\int_0^\infty \frac{\cos(ax)}{x+b} dx = -\text{ci}(ab)\cos(ab) - \text{si}(ab)\sin(ab).$$

**Example 3.4**

We shall evaluate Hecke's integral

$$H(\alpha) = \int_0^\infty \exp\left(-x - \frac{\alpha}{x}\right) \frac{dx}{\sqrt{x}} dx, \ \alpha > 0$$

by using a differential equation (cf. [7]). Differentiation yields

$$H'(\alpha) = -\int_0^\infty \exp\left(-x - \frac{\alpha}{x}\right) \frac{dx}{x\sqrt{x}} dx$$

and then the substitution $x = \alpha/t$ leads to the separable differential equation

$$H'(\alpha) = \frac{-1}{\sqrt{\alpha}} H(\alpha), \ \text{i.e.} \ \frac{dH}{H} = \frac{-d\alpha}{\sqrt{\alpha}}$$

with solution

$$H(\alpha) = M \exp(-2\sqrt{\alpha}), \ M > 0 \ \text{a constant}.$$

Setting here $\alpha = 0$ and using the fact that $H(0) = \Gamma(1/2) = \sqrt{\pi}$ we find

(3.5) $$H(\alpha) = \sqrt{\pi} \exp\left(-2\sqrt{\alpha}\right).$$

**Example 3.5**

Consider the two integrals ([5], p.731)

(3.6) $$U(\alpha) = \int_0^\infty \exp(-x^2)\cos\left(\frac{\alpha^2}{x^2}\right) dx \ \text{and} \ V(\alpha) = \int_0^\infty \exp(-x^2)\sin\left(\frac{\alpha^2}{x^2}\right) dx.$$

We differentiate $U(\alpha)$ and set $y = \alpha/x$ to find

$$U'(\alpha) = -2\int_0^\infty \exp(-x^2)\sin\left(\frac{\alpha^2}{x^2}\right)\frac{\alpha}{x^2}dx = 2\int_\infty^0 \exp\left(\frac{-\alpha^2}{y^2}\right)\sin(y^2)dy$$

$$= -2\int_0^\infty \exp\left(\frac{-\alpha^2}{y^2}\right)\sin(y^2)dy \ .$$

A second differentiation gives

$$U''(\alpha) = -2\int_0^\infty \exp\left(\frac{-\alpha^2}{y^2}\right)\sin(y^2)\frac{-2\alpha}{y^2}dy = 4\int_0^\infty \exp(-x^2)\sin\left(\frac{\alpha^2}{x^2}\right)dx \ ,$$

that is

$$U''(\alpha) = 4V(\alpha).$$

In the same way we compute $V''(\alpha) = -4U(\alpha)$. We define now the complex function $W(\alpha) = U(\alpha) + iV(\alpha)$. This function satisfies the second order differential equation.

$$W'' = -4iW$$

with characteristic equation $r^2 + 4i = 0$ and roots $r_1 = -\sqrt{2} + i\sqrt{2}$ and $r_2 = \sqrt{2} - i\sqrt{2}$. From these roots we construct the general solution to the differential equation

$$W(\alpha) = A\exp(r_1\alpha) + B\exp(r_2\alpha)$$

with parameters $A$ and $B$. Explicitly,

$$W(\alpha) = A\exp(-\sqrt{2}\,\alpha)(\cos\sqrt{2}\,\alpha + i\sin\sqrt{2}\,\alpha) + B\exp(\sqrt{2}\,\alpha)(\cos\sqrt{2}\,\alpha - i\sin\sqrt{2}\,\alpha).$$

At this point we conclude that $B = 0$, since $W(\alpha)$ is a bounded function. Setting $\alpha = 0$ we find $W(0) = A$. At the same time by the definition (3.6) of the above integrals

$$W(0) = U(0) = \int_0^\infty \exp(-x^2)dx = \frac{\sqrt{\pi}}{2}$$

and

$$W(\alpha) = \frac{\sqrt{\pi}}{2}\exp(-\alpha\sqrt{2})(\cos(\alpha\sqrt{2}) + i\sin(\alpha\sqrt{2})).$$

From here, comparing real and imaginary parts we conclude that

(3.7) $$U(\alpha) = \frac{\sqrt{\pi}}{2}\exp(-\alpha\sqrt{2})\cos(\alpha\sqrt{2}), \quad V(\alpha) = \frac{\sqrt{\pi}}{2}\exp(-\alpha\sqrt{2})\sin(\alpha\sqrt{2}).$$

## 4. Advanced techniques

In certain cases we can use the Leibniz Integral Rule [5], [9]:

$$\frac{d}{d\alpha}\int_{\varphi(\alpha)}^{\psi(\alpha)} f(\alpha,x)dx = \int_{\varphi(\alpha)}^{\psi(\alpha)} \frac{d}{d\alpha} f(\alpha,x)dx + f(\alpha,\psi(\alpha))\psi'(\alpha) - f(\alpha,\varphi(\alpha))\varphi'(\alpha),$$

where $f(\alpha,x)$, $\varphi(\alpha)$, $\psi(\alpha)$ are appropriate functions.

**Example 4.1**

We shall evaluate the integral

$$\int_0^1 \frac{\arctan x}{\sqrt{1-x^2}} dx$$

by using the function

$$J(\alpha) = \int_{\varphi(\alpha)}^1 \frac{\arctan(\alpha x)}{\sqrt{1-x^2}} dx$$

where $\alpha > 1$ and

$$\varphi(\alpha) = \sqrt{1 - \frac{1}{\alpha^2}} = \frac{\sqrt{\alpha^2-1}}{\alpha} \quad \text{with} \quad \varphi'(\alpha) = \frac{1}{\alpha^2\sqrt{\alpha^2-1}}.$$

Applying the Leibniz rule we find

$$J'(\alpha) = \int\limits_{\varphi(\alpha)}^{1} \frac{x}{(1+\alpha^2 x^2)\sqrt{1-x^2}} dx - \frac{\arctan\sqrt{\alpha^2-1}}{\alpha\sqrt{\alpha^2-1}}.$$

Let us call this integral $A(\alpha)$. We shall evaluate it by the substitution $1-x^2 = u^2, u > 0$ :

$$A(\alpha) = \int\limits_{\varphi(\alpha)}^{1} \frac{x}{(1+\alpha^2 x^2)\sqrt{1-x^2}} dx = \int\limits_{0}^{\frac{1}{\alpha}} \frac{du}{\alpha^2 + 1 - \alpha^2 u^2}$$

$$= \frac{1}{2\alpha\sqrt{\alpha^2+1}} \ln \frac{\sqrt{\alpha^2+1}+\alpha u}{\sqrt{\alpha^2+1}-\alpha u} \bigg|_{0}^{\frac{1}{\alpha}} = \frac{1}{2\alpha\sqrt{\alpha^2+1}} \ln \frac{\sqrt{\alpha^2+1}+1}{\sqrt{\alpha^2+1}-1}.$$

This function is easy to integrate, as

$$\frac{d}{d\alpha} \ln \frac{\sqrt{\alpha^2+1}+1}{\sqrt{\alpha^2+1}-1} = \frac{-2}{\alpha\sqrt{\alpha^2+1}}$$

and therefore, one antiderivative is

$$\int A(\alpha) d\alpha = \frac{-1}{8} \left( \ln \frac{\sqrt{\alpha^2+1}+1}{\sqrt{\alpha^2+1}-1} \right)^2 .$$

We also have

$$\frac{d}{d\alpha} \arctan\sqrt{\alpha^2-1} = \frac{1}{\alpha\sqrt{\alpha^2-1}}$$

and therefore,

$$\int \frac{\arctan\sqrt{\alpha^2-1}}{\alpha\sqrt{\alpha^2-1}} d\alpha = \frac{1}{2}\left(\arctan\sqrt{\alpha^2-1}\right)^2.$$

Now we can integrate $J'(\alpha)$ to obtain

$$J(\alpha) = \frac{-1}{8}\left(\ln\frac{\sqrt{\alpha^2+1}+1}{\sqrt{\alpha^2+1}-1}\right)^2 - \frac{1}{2}\left(\arctan\sqrt{\alpha^2-1}\right)^2 + C.$$

Using the limit $\lim\limits_{\alpha \to \infty} J(\alpha) = 0$ we find $C = \dfrac{\pi^2}{8}$. Finally,

(4.1) $$\int_{\varphi(\alpha)}^{1} \frac{\arctan(\alpha x)}{\sqrt{1-x^2}} dx = \frac{-1}{8}\left(\ln\frac{\sqrt{\alpha^2+1}+1}{\sqrt{\alpha^2+1}-1}\right)^2 - \frac{1}{2}\left(\arctan\sqrt{\alpha^2-1}\right)^2 + \frac{\pi^2}{8}$$

Setting here $\alpha \to 1$ we find after simple computation

(4.2) $$\int_0^1 \frac{\arctan x}{\sqrt{1-x^2}} dx = \frac{\pi^2}{8} - \frac{1}{2}\left(\ln(1+\sqrt{2})\right)^2.$$

With $\alpha = 2$ in (4.1) we find also

$$\int_{\sqrt{3}/2}^{1} \frac{\arctan(2x)}{\sqrt{1-x^2}} dx = \frac{-1}{8}\left(\ln\frac{\sqrt{5}+1}{\sqrt{5}-1}\right)^2 - \frac{1}{2}\left(\arctan\sqrt{3}\right)^2 + \frac{\pi^2}{8}.$$

**Remark**. Integrating by parts in (4.2) we find

$$\int_0^1 \frac{\arctan x}{\sqrt{1-x^2}} dx = \frac{\pi^2}{8} - \int_0^1 \frac{\arcsin x}{1+x^2} dx$$

and therefore,

$$\int_0^1 \frac{\arcsin x}{1+x^2} dx = \frac{1}{2}\left(\ln(1+\sqrt{2})\right)^2.$$

Using the identity

$$\arctan x = \frac{\pi}{2} - \arctan\frac{1}{x}$$

for $x > 0$, we can write

$$\int_0^1 \frac{\arctan x}{\sqrt{1-x^2}} dx = \frac{\pi^2}{4} - \int_0^1 \frac{\arctan \frac{1}{x}}{\sqrt{1-x^2}} dx.$$

In the second integral we make the substitution $x = 1/t$ to find also

(4.3) $$\int_1^\infty \frac{\arctan t}{t\sqrt{t^2-1}} dt = \frac{\pi^2}{8} + \frac{1}{2}\left(\ln(1+\sqrt{2})\right)^2.$$

This integral was evaluated in [4] independently of (4.2) by using the same method.

**Example 4.2**

We shall evaluate in explicit form the function

$$F(\alpha) = \int_{1/\alpha}^\infty \frac{\ln(\alpha x + \sqrt{\alpha^2 x^2 - 1})}{x(1+x^2)} dx$$

where $\alpha > 0$. Differentiating by the Leibniz rule we compute

$$F'(\alpha) = \int_{1/\alpha}^\infty \frac{dx}{(1+x^2)\sqrt{\alpha^2 x^2 - 1}}$$

(notice that the function $\ln(\alpha x + \sqrt{\alpha^2 x^2 - 1})$ becomes zero for $x = 1/\alpha$).

To solve the integral we first write it in the form

$$F'(\alpha) = \frac{-1}{2} \int_{1/\alpha}^\infty \frac{dx^{-2}}{(x^{-2}+1)\sqrt{\alpha^2 - x^{-2}}}$$

and then we make the substitution $\alpha^2 - x^{-2} = t^2, t > 0$ to get

$$F'(\alpha) = \int_0^\alpha \frac{dt}{1+\alpha^2 - t^2} = \frac{1}{2\sqrt{1+\alpha^2}} \ln \frac{\sqrt{1+\alpha^2} + \alpha}{\sqrt{1+\alpha^2} - \alpha}.$$

This function is easy to integrate, as

$$\frac{d}{d\alpha}\ln\frac{\sqrt{1+\alpha^2}+\alpha}{\sqrt{1+\alpha^2}-\alpha}=\frac{2}{\sqrt{1+\alpha^2}}.$$

Thus we find

$$F(\alpha)=\frac{1}{8}\ln^2\frac{\sqrt{1+\alpha^2}+\alpha}{\sqrt{1+\alpha^2}-\alpha}=\frac{1}{8}\left(\ln(\sqrt{1+\alpha^2}+\alpha)-\ln(\sqrt{1+\alpha^2}-\alpha)\right)^2$$

(the constant of integration is zero, since $F(\alpha)\to 0$ for $\alpha\to 0$). Simplifying this we get

$$F(\alpha)=\frac{1}{2}\ln^2(\sqrt{1+\alpha^2}+\alpha)$$

since $\ln(\sqrt{1+\alpha^2}-\alpha)=\ln\dfrac{1}{\sqrt{1+\alpha^2}+\alpha}=-\ln(\sqrt{1+\alpha^2}+\alpha)$.

Therefore, for any $\alpha>0$,

(4.4) $$\int_{1/\alpha}^{\infty}\frac{\ln(\alpha x+\sqrt{\alpha^2 x^2-1})}{x(1+x^2)}dx=\frac{1}{2}\ln^2(\sqrt{1+\alpha^2}+\alpha).$$

In particular, for $\alpha=1$,

(4.5) $$\int_{1}^{\infty}\frac{\ln(x+\sqrt{x^2-1})}{x(1+x^2)}dx=\frac{1}{2}\ln^2(\sqrt{2}+1).$$

For $\alpha=1/2$ in (4.4) we find

(4.6) $$\int_{2}^{\infty}\frac{\ln(x+\sqrt{x^2-4})}{x(1+x^2)}dx=\ln 2\ln\frac{\sqrt{5}}{2}+\frac{1}{2}\ln^2\frac{\sqrt{5}+1}{2}.$$

**Remark**. Note that the similar integral

(4.7) $$\int_{1}^{\infty}\frac{\ln(x+\sqrt{x^2-1})}{x\sqrt{x^2-1}}dx$$

cannot be evaluated this way. The value of this integral is $2G$, where

$$G = \sum_{n=0}^{\infty} \frac{(-1)^n}{(2n+1)^2} = 1 - \frac{1}{3^2} + \frac{1}{5^2} + \ldots \ ,$$

is Catalan's constant. The substitution $x = \cosh t$ with $\ln(x+\sqrt{x^2-1}) = t$ turns (4,7) into

$$\int_0^{\infty} \frac{t}{\cosh t} dt = 2G$$

which is a well-known result.

# 5. Some theorems

**Theorem A** Suppose the function $f(\alpha, x)$ is defined and continuous on the rectangle $[a,b] \times [c,d]$ together with its partial derivative $f_\alpha(\alpha, x)$. Then

$$\frac{d}{d\alpha} \int_c^d f(\alpha, x) dx = \int_c^d f_\alpha(\alpha, x) dx .$$

In order to apply this theorem in the case of improper integrals we have to require uniform convergence of the integral with respect to the variable $\alpha$. A simple sufficient condition for uniform convergence is the following theorem

**Theorem B** Suppose $f(\alpha, x)$ is continuous on $[a,b] \times [0, \infty)$ and $g(x)$ is integrable on $[0, \infty)$. If

$$|f(\alpha, x)| \leq g(x)$$

for all $a \leq \alpha \leq b$ and all $x \geq 0$, then the integral

$$\int_0^{\infty} f(\alpha, x) dx$$

is uniformly convergent on $[a,b]$.

**Theorem C** Suppose the function $f(\alpha, x)$ is continuous on $[a,b] \times [0,\infty)$ together with its partial derivative $f_\alpha(\alpha, x)$. In this case

$$\frac{d}{d\alpha} \int_0^\infty f(\alpha, x) dx = \int_0^\infty f_\alpha(\alpha, x) dx$$

when the first integral is convergent and the second is uniformly convergent on $[a,b]$.

The case of improper integrals on finite intervals is treated in the same way. For details and proofs we refer to [5], [7], [9], and [12]. The book [5] presents the Leibniz rule in full detail.

# References

[1]  **Khristo N. Boyadzhiev**, Some integrals related to the Basel problem, *SCIENTIA. Series A: Mathematical Sciences,* 26 (2015), 1-13. Also arXiv:1611.03571 [math.NT]

[2]  **Khristo N. Boyadzhiev,** On a series of Furdui and Qin and some related integrals, 2012, arXiv:1203.4618v3 [math.NT]

[3]  **Khristo Boyadzhiev, Hans Kappus,** Solution to problem E 3140, *Amer. Math. Monthly*, 95 (1), (1988), 57-59.

[4]  **Hongwey Chen**, Parametric differentiation and integration, *Internat. J. Math. Ed. Sci. Tech*. 40 (4) (2009), 559-579.

[5]  **G. M. Fikhtengol'ts**, *A Course of Differential and Integral Calculus* (Russian), Vol. 2, Nauka, Moscow, 1966.

   Abbreviated version in English: *The Fundamentals of Mathematical analysis*, Vol. 2, Pergamon Press, 1965.

[6]  **I. S. Gradshteyn and I. M. Ryzhik**, *Tables of Integrals, Series, and Products*, Academic Press, 1980.

[7]  **Omar Hijab**, *Introduction to Calculus and Classical Analysis*, Springer, 1997.

[8]  **A. P. Prudnikov, Yu. A. Brychkov, O. I. Marichev**, *Integrals and Series, vol.1 Elementary Functions*, Gordon and Breach 1986.

[9]  **Ioannis Markos Roussos**, *Improper Riemann Integrals*, Chapman and Hall/CRC, 2014.

[10] **Joseph Wiener**, Differentiation with respect to a parameter, *College Mathematics Journal,* 32, No. 3.(2001), pp. 180-184.

[11] **66th Annual William Lowell Putnam Mathematical Competition**, *Math. Magazine*, 79 (2006), 76-79.

[12] **Frederick. S. Woods**, *Advanced Calculus: A Course Arranged with Special Reference to the Needs of Students of Applied Mathematics*. New Edition, Ginn and Co., Boston, MA, 1934.

[13] **Aurel J. Zajta, Sudhir K. Goel**, Parametric integration techniques, *Math. Magazine*, 62(5) (1989), 318-322.